\documentclass[lang = american]{ems-icm} 

\theoremstyle{plain}
\newtheorem{thm}{Theorem}[section]
\newtheorem{lemma}[thm]{Lemma}
\newtheorem{prop}[thm]{Proposition}
\newtheorem{cor}[thm]{Corollary}
\newtheorem{conjec}[thm]{Conjecture}
\newtheorem{qtn}[thm]{Question}
\newtheorem{prob}[thm]{Problem}

\theoremstyle{definition}
\newtheorem{defn}[thm]{Definition}

\newtheorem{rem}[thm]{Remark}

\DeclareMathOperator{\SL}{SL}

 \DeclareMathOperator{\SO}{SO}

\DeclareMathOperator{\SU}{SU} 
\DeclareMathOperator{\SP}{Sp}

\DeclareMathOperator{\Isom}{Isom}

\DeclareMathOperator{\Diff}{Diff}

\newcommand{\Ga}{\Gamma}

\newcommand{\Gam}{\Gamma}

\newcommand{\bfH}{\mathbf{H}}

\newcommand{\Hy}{\ensuremath{\mathbb{H}}}
\newcommand{\bfP}{\ensuremath{\mathbb{P}}}
\newcommand{\CH}{\ensuremath{\mathbb{CH}}}

\newcommand{\R}{\ensuremath{\mathbb{R}}}
\newcommand{\Z}{\ensuremath{\mathbb{Z}}}
\newcommand{\Ca}{\ensuremath{\mathbb{C}}}

\newcommand{\Ta}{\ensuremath{\mathbb{T}}}
\newcommand{\bfJ}{\mathbf{J}}
\newcommand{\bfL}{\mathbf{L}}

\numberwithin{equation}{section}

\begin{document}

\volumetitle{ICM 2022} 

\title{Rigidity, lattices and invariant measures beyond homogeneous dynamics}
\titlemark{Rigidity, lattices and measures}

\emsauthor{1}{David Fisher}{D.~Fisher}

\emsaffil{1}{Department of Mathematics, Indiana University, Bloomington IN 47401 USA  \email{fisherdm@indiana.edu}}


\begin{abstract}
This article discusses two recent works by the author, one with Brown and Hurtado on Zimmer's conjecture and one with Bader, Miller and Stover
on totally geodesic submanifolds of real and complex hyperbolic manifolds.  The main purpose of juxtaposing these two very disparate sets
of results in one article is to emphasize a common aspect:  that the study of invariant and partially invariant measures outside the homogeneous setting is important
to questions about rigidity in geometry and dynamics.  I will also discuss some open questions including some that seem particularly compelling in light of this
juxtaposition.
\end{abstract}

\maketitle

\section{Introduction}

This article focuses on some recent developments concerning the rigidity of discrete subgroups of Lie groups.  This area has long had deep connections
with ergodic theory and dynamical systems going back to seminal work of Furstenberg, Margulis, Mostow and Zimmer \cite{MR0284569, MR0422499, MargulisICM ,Mostow-book, ZimOE}.   Here we emphasize the role of invariant measures
for groups and their subgroups.  From one point of view, a key step in Margulis' proof of his superrigidity theorem can be written as finding an invariant measure for
a subgroup in a group action. Invariant measures have also long been a key object of study in homogeneous dynamics.   Rigidity of discrete groups and homogeneous dynamics have long been allied and interacting fields, but the developments here are part of a strengthening of those connections, particularly in terms of proving rigidity results directly through the study of invariant measures in inhomogeneous settings, using techniques, ideas and results from homogeneous dynamics.  In particular, both results study dynamical systems with homogeneous factors and the dynamics on the homogeneous factor help control the dynamics on the total system.

Let $H$ be a Lie group, $\Lambda <H$ a discrete subgroup and $S<H$ a subgroup. Homogeneous dynamics is most often the study of invariant measures, orbit closures and equidistribution for the $S$ action on $H/\Lambda$.  This area has been quite fruitful with applications to areas as diverse as number theory, geometry and physics.  Here we
are more often concerned with dynamical systems where the space $H/\Lambda$ is replaced by a space that is not homogeneous.  Perhaps the most famous example of this is the action of $SL(2,\mathbb{R})$ on the moduli space of quadratic differentials on a surface.  This area is not our topic here, but the fruitful importation of ideas from homogeneous dynamics to this area has a long history, starting with work of Veech and currently culminating in the work of Eskin, Mirzakhani and Mohammadi \cite{Veech, EM, EMM}.  Some ideas arising in this setting have been pushed even further into the inhomogeneous world by work of Brown-Rodriguez Hertz and ongoing work of Brown-Eskin-Filip \cite{BRH}.   More closely allied with the developments here is work of Katok, Kalinin and Rodriguez Hertz on invariant measures for actions of higher rank abelian groups on compact manifolds \cite{KKRH}.

I want to point to one other principle that has played a key role in many developments, which is the notion of stiffness, first formalized by Furstenberg \cite{FurstenbergStiff}.  When studying acitons of amenable groups, it suffices to consider invariant measures.  When the acting group $S$ is not amenable, then to understand the dynamics it is important to study the broader class of stationary measures.  Furstenberg called an action of a group $S$ on a space $X$  \emph{stiff} if every stationary measure for $S$ was in fact invariant.  Even before Furstenberg defined the term, Nevo and Zimmer had considered the case where $S$ is a higher rank simple Lie group acting on an arbitrary measure space $X$ and given criteria for stiffness in terms of measurable projective quotients \cite{NevoZimmer}. More recently, Benoist and Quint proved dramatic results on stiffness in the homogeneous setting that were inspirational for the work of Eskin-Mirzakhani mentioned above \cite{BenoistQuint}.

In this article we point to directions where homogeneous dynamics techniques are applied outside the homogeneous setting to prove rigidity results about discrete groups and their actions. And in particular questions where we need to move beyond the question of stiffness.  The first work I will describe, joint with Brown and Hurtado, resolves important cases of Zimmer's conjecture.  A key ingredient in our work is to move beyond stiffness and consider an even larger class of measures than the stationary ones. In this context, the stationary measures can be made to correspond to the invariant ones for some subgroup $P<S$ but our proof requires understanding something about the invariant measures for a much smaller subgroup $A<P$.   The other results I will focus on concern totally geodesic submanifolds of real and complex hyperbolic manifolds of finite volume and are joint works with Bader, Lafont, Miller and Stover.  In this setting, once again a key object is to construct certain measures that are invariant under subgroups of the acting group.  In this work, the class of measures studied is simply different than the stationary ones, not more or less general.

Overall, I think the results mentioned above, as well as numoreous results by other authors, point to the development of a broad area of research in which the study of invariant measures beyond homogeneous dynamics has broad implications for rigidity theory.  And that many of these developments will spring from broadening the efficacy of ideas that originate in homogeneous dynamics.  It feels too early to attempt a survey of these developments, so I will restrict myself to an account of these developments in my own work.  Along the way I will point to open questions inspired by that work, some of which fits this introductory framework and some of which do not.

\section{Zimmer's conjecture and the Zimmer program}

In this section I will discuss some aspects of recent work with Brown and Hurtado on Zimmer's conjecture and also discuss the implications for further work.  In the course of this, I will point to some work in progress with Melnick concerning examples and also point to an old paper of Uchida which has important implications for the Zimmer program and which seems too little known.  For a different take on some results and questions discussed here, see Brown's contribution in these proceedings \cite{BrownICM}.

\subsection{Zimmer's conjecture}

Throughout this subsection $G$ will be a simple Lie group with real rank at least $2$ and $\Gamma <G$ will be a lattice.  The reader will lose little by considering
the case of $G=\SL(n,\R)$ and $\Gamma = \SL(n,\Z)$ with $n>2$.  In \cite{ZimmerICM,ZimmerBulletin}, Zimmer laid out a program for understanding $\Gamma$ actions on a compact manifold $M$ preserving a volume form $\omega$.  The base case of this program is to show that for any homomorphism $\rho:\Gamma \rightarrow \Diff(M,\omega)$, the image $\rho(\Gamma)$ always preserves a smooth Riemannian metric if $\dim(M)$ is less than the dimension of the minimal real $G$ representation.  This conjecture was motivated by Zimmer's cocycle superrigidity theorem, which in this context produced a measurable Riemannian metric whose associated volume form was $\omega$.

Already in the 1990's, this conjecture about low dimensional actions had been transported out of the volume preserving setting by numerous works, particularly concerning actions on the circle, see e.g. \cite{DH, Ghys}.  The work with Brown and Hurtado proved this more general conjecture in many cases, here I state only the following special case.

\begin{thm}
\label{thm:main}[Brown, Fisher, Hurtado]
Let $\Gamma$ be a lattice in $\SL(n,\R)$, let $M$ be
a compact manifold and let $\rho:\Ga \rightarrow \Diff(M)$ be a homomorphism.
Then
\begin{enumerate}
  \item if $\dim(M) < n-1$, the image of $\rho$ is finite;
  \item if $\dim(M) < n$ and $\rho(\Gamma)$ preserves a volume form on $M$,
  then the image of $\rho$ is finite.
\end{enumerate}
\end{thm}

While I include a very rough sketch of ideas in the proof, more detailed outlines can be found in \cite{BrownBook, FisherZimUpdate, FisherNotices}.  A key step is showing that
$\Gamma$ acts with \emph{subexponential growth of derivatives}.  We let $l$ be any word length on the group $\Gamma$.

\begin{defn}
\label{def:subexp}
Let $\rho: \Gamma \rightarrow \Diff(M)$ be an action of a finitely generated group $\Gamma$ on a compact manifold $M$.
We say $\rho$ has {\em subexponential growth of derivatives} if for every $\epsilon >0$ there exists $C>0$ such that
$$\sup_{x \in M} \|D\rho(\gamma)_x\| < C e^{\epsilon l(\gamma)}.$$
\end{defn}

The idea that controlling growth of derivatives is pivotal was long known, see the discussion in \cite[End of Section 7]{FisherZimUpdate}.  However, a particular novelty introduced in \cite{BFH1} is that we use the strong property $(T)$ of Lafforgue to convert subexponential growth of derivatives to an invariant Riemannian metric without any further hypothesis.  For much more on strong property $(T)$, see de la Salle's contribution in these proceedings \cite{delaSalleICM}.

The approach to proving subexponential growth of derivatives in \cite{BFH1, BFH2, BFH3} was also distinct from previous ideas on Zimmer's conjecture but did have some classical inspirations as well as a more closely related one in \cite{HurtadoBurnside}.  We give a somewhat ad hoc definition of zero Lyaponov exponents for a group action here, other definitions are possible and most are somewhat weaker than this, but this one suffices for current purposes.

\begin{defn}
\label{defn:zeroexp}
Let $\rho: \Gamma \rightarrow \Diff(M)$ be an action of a finitely generated group on a compact manifold.  Let $\mu$ be a measure on $M$.  We say $\rho$ has {\em zero first Lyapunov exponent} for $\mu$ if

$$\lim_{l(\gamma) \rightarrow \infty} \frac{\ln \| D\rho(\gamma)_x\|}{l(\gamma)} = 0$$

\noindent for $\mu$ almost every $x$ in $M$.
\end{defn}

Essentially, the proof of Theorem \ref{thm:main} consists of showing that exponential growth of derivatives must be witnessed by a positive Lyapunov exponent for some $\Gamma$ invariant measure and then seeing that this contradicts Zimmer's cocycle superrigidity theorem.   The motivation is probably most easily encapsulated by this classical propostion.

\begin{prop}
\label{prop:motivation}
Let $M$ be a compact manifold and let $\rho: \Z \rightarrow \Diff(M)$ be the action generated by a single diffeomorphism $f$.
Then $\rho$ has subexponential growth of derivatives if and only if $\rho$ has zero first Lyapunov exponent for every $f$ invariant
measure $\mu$.
\end{prop}

\noindent  The proposition is not easily adapted to group actions partly because if the group is not amenable, there may be no invariant measures at all.
One attempt to remedy this would be to consider stationary measures and random Lyapunov exponents, but there is no useful analogue of the proposition in that context.
The issue that arises is that the random Lyaponuv exponent might vanish while exponential growth of derivatives still occurs along a very thin set of trajectories
in the acting group.  There is a  weaker statement that is true and even useful in some contexts that is implicit in \cite[Section 3]{HurtadoBurnside}.  That argument does produce a measure on a skew product over a shift space.  However, the measure produced when projected to the shift space is not an independent and identically distributed but quite arbitrary and does not fall into the usual context of random dynamics, stationary measures and stiffness.

For the work on Zimmer's conjecture, we take a long detour which I will only describe a small part of in the next subsection, to make clear the connections to homogeneous dynamics and also to make clear why it is not enough to prove stiffness of the action.

\subsection{Measures in the proof of Zimmer's conjecture}
\label{subsection:zimproof}

The first step in the proof of Theorem \ref{thm:main} uses the notion of an induced action, a variant of induced representations due to Mackey. This notion is also similar to the construction of flat bundles.  If $\Gamma$ acts on a manifold $M$ via a homomorphism $\rho: \Gamma \rightarrow \Diff(M)$, then we can build a $G$ action on a manifold $(G \times M)/\Gamma$.  This can be specified just by specifying commuting $G$ and $\Gamma$ actions on $G \times M$.  We do this by the formula $$g(g_0, m) \gamma = (gg_0\gamma^{-1}, \rho(\gamma)m).$$
\noindent Note that there is a $G$ equivariant map $\pi: (G\times M)/\Gamma \rightarrow G/\Gamma$ and this map exhibits $(G \times M)/\Gamma$ as a fiber bundle over $G/\Gamma$ with fiber $M$.  Also note that the tangent bundle to $(G \times M)/\Gamma$ admits a $G$ invariant subbundle consisting of directions tangent to fibers of the projection $\pi$, i.e.  $(G \times TM)/\Gamma \subset T(G\times M)/\Gamma$.

We then define {\em fiberwise zero first Lyapunov exponent} and {\em fiberwise subexponential growth of derivatives} for the $G$
action on $(G \times M)/\Gamma$ by restricting all derivatives in the definitions to the invariant subbundle $(G \times TM)/\Gamma \subset T(G\times M)/\Gamma$.
It is a relatively easy exercise to see that if $\Gamma < G$ is cocompact, then subexponential growth of derivatives for the $\Gamma$ action is equivalent to fiberwise subexponential growth of derivatives for the $G$ action on $(G \times M)/\Gamma$.  The situation when $G/\Gamma$ has finite volume but is not compact is considerably more complicated and we do not discuss it here.  At this point, the structure of Lie groups begins to play an important role.  It turns out that $G$ can always be written as a product $KAK$ where $K$ is compact and $A$ is abelian. For $\SL(n,\R)$, these groups are $K=SO(n)$ and $A$ the group of diagonal matrices of determinant one.  Since $K$ is compact, we can average
any Riemannian metric on $(G \times M)/\Gamma$ over the $K$ action and obtain a $K$ invariant metric.   This means that for the action of $G$ any growth of derivatives that we see comes entirely from the action of $A$.  Modifying the proof of Proposition \ref{prop:motivation} and retaining the notation and terminology above, we prove

\begin{lemma}
\label{lemma:aok}
Given $\rho: \Gamma \rightarrow \Diff(M)$ then either $\rho$ has subexponential growth of derivatives or there is an $A$ invariant
measure $\mu$ on $(G \times M)/\Gamma$ with non-zero fiberwise first Lyapunov exponent for some element $a$ in $A$.
\end{lemma}

If $\mu$ were in fact $G$ invariant, then this can be seen to contradict Zimmer's cocycle superrigidity theorem.   A key point that
makes it possible to prove Lemma \ref{lemma:aok} is that $A$ is abelian and so amenable.

We proceed by proving that $\mu$ can be replaced by a measure that is in fact $G$ invariant.  This is done in two steps. First we average the measure over certain subgroups of $G$ to produce a measure $\mu'$ whose projection $\pi_*\mu'$ to $G/\Gamma$ is Haar measure. The difficulty here is to do the averaging while retaining that $\mu'$ is $A$ invariant and that some $a$ in $A$ has positive first Lyapunov exponent for $\mu'$.   After this step, we can use a result of Brown, Rodriguez Hertz and Wang together with some algebraic computations to show that $\mu'$ is in fact $G$ invariant \cite{BRHW}.  This contradiction shows that $\rho$ does in fact have subexponential growth of derivatives.

The step of averaging the measure $\mu$ to produce $\mu'$ makes extensive use of homogeneous dynamics and, in particular, work
of Ratner and Shah \cite{MR1135878, MR1106945, MR1291701}.  The work of Brown, Rodriguez Hertz and Wang pivots on relations between invariant measures and entropy and in particular on an extension of the important work of Ledrappier and Young \cite{MR819556}.  A key ingredient in both parts is the theory of Lyapunov exponents and particularly the fact that for actions of an abelian group $A$, Lyapunov exponents give rise to linear functionals on $A$.

In this context, stationary measures can be thought of roughly as just measures invariant under $P$ and not all of $G$.  It turns out not to be possible to start a proof
by considering only $P$ invariant measures, instead of the wider class of $A$ invariant measures.  A priori, exponential growth of derivatives is only witnessed on
some sequence $g_n$ of elements in $G$ and a naive rewriting of the proof of Proposition \ref{prop:motivation} only gives an exponent along that particular sequence. One can choose this sequence to be in $P$ by essentially the same reasoning by which we choose it to be in $A$.  But at that point, it is then very hard to proceed since Lyapunov exponents on $P$ do not a priori have any structure analogous to their structure as linear functionals on $A$.  A main observation is that one can in fact choose $g_n$ to be in $A$ and use this to produce Lyapunov exponent first for a measure invariant under a $1$-parameter subgroup of $A$ and then by averaging under all of $A$.  While our argument has been rewritten by An, Brown and Zhang as then producing a $P$ invariant measure from this $A$ invariant measure, at this step one has to argue using a great deal of homogeneous dynamics and using a rather intricate averaging procedure \cite{ABZ}.

\subsection{Other results, future directions}

The joint work with Brown and Hurtado classifies actions of lattices in $\SL(n,\R)$ in dimension at most $n-2$ and volume preserving actions in dimension $n-1$.
A more recent result of Brown, Rodriguez Hertz and Wang completes the picture through dimension $n-1$, showing that in general an action of a lattice in $\SL(n, \R)$
on an $n-1$ dimensional manifold either factors through a finite group or extends to an $\SL(n,\R)$ action.   It is easy to see that there are exactly two actions of $\SL(n,\R)$ on $n-1$ manifolds, namely the action on $\bfP(\R^n)$ and lift of that action to $S^{n-1}$. Clearly a natural question would be to also classify actions in dimension $n$, but
this becomes surprisingly harder, mainly because there are many more examples.

Already actions of $\SL(n,\R)$ on $n$ manifolds are quite complicated and not fully classified.  A remarkable and little known paper of Uchida classifies analytic actions of $\SL(n,\R)$ on $S^n$ and the parameter space turns out to be infinite dimensional.  In work in progress with Melnick, we are extending this to a classification of all analytic actions of $\SL(n,\R)$ on $n$ manifolds and may also produce a smooth classification, though there are missing ingredients at the moment.  For lattice actions, there are more examples known, constructed first by Katok and Lewis \cite{KatokLewis} and studied by the author with various coauthors \cite{BenvenisteFisher, FisherWhyte}.  In addition, some ideas from the work of Uchida allow us to adapt a continuous construction described by Farb and Shalen and show that it can be done analytically \cite{FS, Uchida}.

The current conjectural picture of actions of $\SL(n,\R)$ lattices on $n$ manifolds that results from these developments is quite complicated and we will not attempt to describe it here.  Instead we state a conjecture in dimension $n$ about stiff actions.

\begin{conjec}
\label{stiffdimn}
Let $\Gam < \SL(n,\R)$ be a lattice and assume $\Gamma$ acts on an $n$ manifold stiffly.  Then either the action factors through a finite group
or the action lifts to a finite cover where it is smoothly conjugate to an affine action of a finite index subgroup of $\SL(n,\Z)$ on $\Ta^n$.
\end{conjec}

For non-stiff actions, one can formulate a conjecture about actions of lattices in $\SL(n,\R)$ on $n$ manifolds but the statement
becomes quite involved.  The pivotal fact that one would need to even begin classifying actions is summarized by

\begin{conjec}
\label{conject:dimnPmeasures}
Let $\Gam <\SL(n,\R)$ be a lattice and assume $\Gamma$ acts on an $n$ manifold $M$.  Assume that the induced action on $(G\times M)/\Gamma$
admits a $P$ invariant measure than is not $G$ invariant.  Then the support of this measure is of the form $(G \times N)/\Gamma$ where $N$
is an embedded $n-1$ sphere or $\bfP(\R^n)$ in $M$.
\end{conjec}

One can even weaken the hypotheses to consider $A$ invariant measures that are not $G$ invariant and it seems that the only additional
possibility will be Haar measure on closed $A$ orbits in $(G\times M)/\Gamma$.
While I don't state the full conjecture here, I do hope to include it in future work.  An added difficulty is that it seems one can obtain analytic actions with an open subset where the action is analytically conjugate to the one that extends to the action
of $\SL(n, \R)$ on $\R^n \backslash \{0\}$.  These subsets do not support any invariant probability measures and new ideas are definitely needed
to capture this behavior in a classification.   Another key step in completing a classification would involve designing an equivariant surgery that cuts along the invariant $N$ from Conjecture \ref{conject:dimnPmeasures} and simplifies the resulting manifold and action. A full classification may be considerably easier if one assumes the action is volume preserving.
This rules out the ``bad'' open sets just described and would also considerably simplify the required surgery operations.

In fact, in the current state of knowledge, one expects that a variant of Conjecture \ref{stiffdimn} might hold much more generally.
For this we require a definition from \cite{FisherMargulisLRC}.

\begin{defn}
\label{definition:affine}
 {\textbf 1} Let $A$ and $D$ be topological
groups, and $B<A$ a closed subgroup. Let
$\rho:D{\times}A/B{\rightarrow}A/B$ be a continuous action. We
call $\rho$ {\em affine}, if, for every $d{\in}D$ there is a
continuous automorphism $L_d$ of $A$ and an element $t_d{\in}A$
such that $\rho(d)[a]=[t_d{\cdot}L_d(a)]$.

\smallskip
\noindent {\textbf 2} Let $A$ and $B$ be as above.  Let $C$ and $D$ be
two commuting groups of affine diffeomorphisms of $A/B$, with $C$
compact. We call the action of $D$ on $C{\backslash}A/B$ a {\em
generalized affine action}.

\smallskip
\noindent {\textbf 3} Let $A$, $B$, $D$ and $\rho$ be as in $1$ above.
Let $M$ be a compact Riemannian manifold and
$\iota:D{\times}A/B{\rightarrow}\Isom(M)$ a $C^1$ cocycle.  We
call the resulting skew product $D$ action on $A/B{\times}M$ a
{\em quasi-affine action}. If $C$ and $D$ are as in $2$, and
$\alpha:D{\times}C{\backslash}A/B{\rightarrow}\Isom(M)$ is a $C^1$
cocycle, then we call the resulting skew product $D$ action on
$C{\backslash}A/B{\times}M$ a {\em generalized quasi-affine
action}.
\end{defn}

We note that generalized quasi-affine actions for $D$ a higher rank simple group
or a lattice in such a group might be more constrained than it first appears.
It seems at least possible that there are considerable restrictions
on the cocycle into $\Isom(M)$ defining the skew product.  Very partial results
in this direction are obtained by Witte Morris and Zimmer in \cite{DaveBob}.

\begin{qtn}
\label{stiffgeneral}
Assume $G$ is a higher rank simple Lie group and $\Gam<G$ is a lattice and $M$ is a compact manifold.  Let $\rho:\Ga \rightarrow \Diff(M)$
be a stiff action.  Is the action generalized quasi-affine?
\end{qtn}

\noindent  A negative answer to Question \ref{stiffgeneral} would require a genuinely new idea.  All existing constructions of actions which
are not generalized quasi-affine involve cutting and pasting along certain singular divisors and these singular divisors always carry
stationary measures that are not invariant.

All of our understanding of these examples suggest the following

\begin{conjec}
Assume the set up of Question \ref{stiffgeneral}. Instead of assuming stiffness, assume that the action preserves a rigid geometric structure or that
a single element admits a dominated splitting, then the action is generalized quasi-affine.
\end{conjec}

For all of the conjectures and questions just mentioned, one might try to start with $P$ invariant measures in the induced action, but the proof
of Zimmer's conjecture does indicate that it is perhaps more fruitful to start by studying $A$ invariant measures.  It also indicates that one
intermediate goal might be producing some uniform hyperbolicity, such as a dominated splitting.  Subexponential growth of derivatives
is exactly the uniform absence of hyperbolicity.  The best evidence for the dynamical form of this conjecture is in dimension $n$ in a recent paper of my student
Homin Lee \cite{Homin}.  There is more plentiful evidence for the geometric form, going back to results of Zimmer, but the question remains mostly open see
\cite[Section 6]{FisherSurvey}.

In addition, one might ask for some analogue of Conjecture \ref{conject:dimnPmeasures}.  The analogue might well have the identical
statement but it is not as clear what $N$ should occur here. Examples first constructed by Benveniste shows that the
cutting and pasting may occur on along much more complicated submanifolds in general, see e.g. \cite{FisherDeformation, BenvenisteFisher}.
In these examples, the $P$ invariant but not $G$ invariant measures are in fact supported on more complicated
sets and not on manifolds of the form $G/Q$ where $Q$ is a parabolic of $G$.  Instead one gets for example sets that
are $G/Q$ bundles over $G/\Gamma$ for some lattice.  One can also build examples where the cutting and pasting
occurs along a high dimensional sphere or projective space and the $P$ invariant measures are supported on
a very low dimensional submanifold of that space simply by doing the Katok-Lewis example for some large
value of $N$ and restricting the action to $\SL(n,\Z)< SL(N,\Z)$ for some $3 \leq n << N$.  In these examples
one will also end up with invariant open subsets where the induced action admits no stationary probability measure.

Because of this complexity, it quickly becomes hard to state a general conjecture in high dimensions precisely.
A formulation favored by Zimmer and later Labourie is that the action is homogeneneous on an open dense set
or perhaps even built of locally homogeneous pieces.   For a different conjecture, concerning ways in which
one might expect the $\Gamma$ action to extend locally to a $G$ action, see \cite[Conjecture 5.6]{FisherMargulissurvey}.

Before ending this section, I will recall that while the results we can prove for Zimmer's conjecture are sharp for lattices $\SL(n,\R)$
and $\SP(2n, \R)$ they are not sharp for lattices in other simple and semi-simple Lie groups.  There are two obstructions that arise,
each of which is serious.  For non-split groups, including even $\SL(n, \Ca)$ and $\SL(n, \Hy)$, there is an issue arising where we employ
the work of Brown, Rodriguez Hertz and Wang \cite{BRHW}.  That work only ``sees" the number of roots of a simple Lie group and not the dimensions
of the root subspaces.  To see the impact on dimensions where we can prove a result, recall that any simple Lie group $G$ contains a maximal
$\R$-split subgroup $G'$ of the same $\R$-rank.  If we let $Q$ be a maximal parabolic of highest dimension $G$ and $Q'$ a maximal parabolic
of highest dimension in $G'$, then one  expects that for all lattices $\Gam <G$ that all smooth actions on compact manifolds are isometric below $\dim(G/Q)$ but our methods only
prove this for $\dim(G'/Q')$.  To see the effect in practice, one should consider say $\SO(m,n)$ for $m<n$ and note that the maximal
$\R$-split subgroup is $\SO(m,m+1)$.  This shows that this gap between expected results and what we can prove can be arbitrarily large.
For $\Ca$-split groups, An, Brown and Zhang have announced a remedy to this issue, but it appears their solution for that case is
not sufficiently robust to overcome the problem in general \cite{ABZ}.

There is another gap that arises from the fact that our proofs seem, in most cases, to be much better in the context
where one does not preserve a volume form.  Our techniques always only manage to constrain volume preserving actions
in one dimension more than they preserve non-volume preserving ones.  This is because all we are using about volume preserving
actions is the single linear condition an invariant volume form imposes on Lyapunov exponents.  Conjecturally all volume
preserving  actions should be isometric below the dimension $d$ of the minimal $G$ representation, while all actions should
only be isometric below the dimension of $G/Q$ where $Q$ is a maximal parabolic of largest dimension.  For the real
split form of $E_8$ these numbers are $248$ and $57$.  Surprisingly this is the largest gap that occurs between the
conjectured minimal dimensional non-isometric  volume preserving action and the conjectured minimal dimensional non-isometric action.

\begin{prob}
Complete the proof of Zimmer's conjecture in general by overcoming the two problems just discussed.
\end{prob}

\noindent If I had to guess, I would say that the second problem is considerably harder than the first to overcome.  There is
no clear robust dynamical behavior to exploit to resolve the problem.

\section{Totally geodesic manifolds and rank one symmetric spaces}

This section concerns recent results by Bader, the author, Miller and Stover, motivated by questions of McMullen and Reid in the case of real hyperbolic manifolds.  Throughout this section a geodesic submanifold will mean a closed immersed, totally geodesic submanifold. (In fact all results can be stated also for orbifolds but we ignore this technicality here.) A geodesic submanifold is {\emph{maximal} if it is not contained in a proper geodesic submanifold of smaller codimension.

For arithmetic manifolds, the presence of one maximal geodesic submanifold can be seen to imply the existence of infinitely many.  The argument involves lifting the submanifold $S$ to a a finite  cover $\tilde M$ where an element $\lambda$ of the commensurator acts as an isometry.  It is easy to check that for most choices of $\lambda$, the submanifold $\lambda(S)$ can be pushed back down to a geodesic submanifold of $M$ that is distinct from $S$.  This was perhaps first made precise in dimension $3$ by Maclachlan--Reid and Reid \cite{MRTG, ReidTG}, who also exhibited the first hyperbolic $3$-manifolds with no totally geodesic surfaces.

In the real hyperbolic setting the main result from \cite{BFMS} is

\begin{thm}[Bader, Fisher, Miller, Stover]
\label{thm:bfms1}
Let $\Gam$ be a lattice in $\SO_0(n, 1)$. If the associated locally symmetric space contains infinitely many maximal geodesic submanifold of dimension at least $2$, then $\Gam$ is arithmetic.
\end{thm}

\begin{rem} \begin{enumerate}

\item The proof of this result involves proving a superrigidity theorem for \emph{certain} representations of the lattice in $\SO_0(n,1)$.  As the conditions required become a bit technical, we refer the interested reader to
    \cite{BFMS}.  The superrigidity is proven using ideas and methods introduced in \cite{BaderFurman1}.

\item At about the same time, Margulis and Mohammadi gave a different proof for the case $n=3$ and $\Gam$ cocompact \cite{MM}.  They also proved  a superrigidity theorem, but both the statement and the proof are quite different than in \cite{BFMS}.

\item A special case of this result was obtained a year earlier by the author, Lafont, Miller and Stover \cite{FLMS}. There we prove that  a large class of non-arithmetic manifolds have only finitely many maximal totally geodesic submanifolds.  This includes all the manifolds constructed by Gromov and Piatetski-Shapiro but not the examples constructed by Agol and Belolipetsky-Thomson.

\end{enumerate}
\end{rem}

 Theorem \ref{thm:bfms1} has a reformulation entirely in terms of homogeneous dynamics and homogenenous dynamics play a key role in the proof.  It is also interesting that a key role is also played by dynamics that are not quite homogeneous but that take place on a projective bundle over the homogeneous space $G/\Gamma$.  In fact, the work can be used to give a classification of invariant measures for certain subgroups $W<SO_0(n,1)$ on these projective bundles.

Even more recently the same authors have extended this result to cover the case of complex hyperbolic manifolds.

\begin{thm}[Bader, Fisher, Miller, Stover] \label{thm:charithmetic}
Let $n \ge 2$ and $\Gam < \SU(n,1)$ be a lattice and $M=\CH^n/\Gamma$. Suppose that $M$ contains infinitely many maximal totally geodesic submanifolds of dimension at least $2$, then $\Gam$ is arithmetic.
\end{thm}

As before this is proven using homogeneous dynamics, dynamics on a projective bundle over $G/\Gamma$, and a superrigidity theorem.  Here the superrigidity theorem is even more complicated than before and depends also on results of Simpson and Pozzetti \cite{Simpson, Pozzetti}.  A very different proof for the case where the totally geodesic submanifolds are all assumed to be complex submanifolds was given very shortly after ours by Baldi and Ullmo \cite{BaldiUllmo}.  There is almost no overlap of ideas between the two proofs, theirs characterizes the totally geodesic submanifolds in terms of special intersections and then studies them using Hodge theory and $o$-minimality.

The results in this section provide new evidence that totally geodesic manifolds play a very special role in non-arithmetic lattices and perhaps provide some evidence that the conventional wisdom on Questions \ref{question:tg} and \ref{qtn:nonarithmeticSU} below should be reconsidered.

\subsection{Other results and open questions}
\label{subsub:constructionsandquestions}

It is important to preface this section by saying that for all semisimple Lie groups $G$ other than $\SO(n,1)$ for $n>2$ and $\SU(n,1)$ for $n>1$, we have an essentially complete classification of lattices in $G$.  For $S0(2,1)=SU(1,1)$, the lattices are exactly the fundamental groups of hyperbolic surfaces of finite volume, which have been understood for
quite some time.  For all the remaining groups, all lattices are arithmetic.  I will discuss the known construction of non-arithmetic lattices in $\SO(n,1)$ and $\SU(n,1)$. To begin slightly out of order, I emphasize one of the most important open problems in the area.

\begin{qtn}
\label{qtn:nonarithmeticSU}
For what values of $n$ does there exist a non-arithmetic lattice in $\SU(n,1)$?
\end{qtn}

The answer is known to include $2$ and $3$.  The first examples were constructed by Mostow in \cite{MostowPJM} using reflection group techniques. The list was slightly expanded by Mostow and Deligne using monodromy of hypergeometric functions \cite{DeligneMostow, MostowGPL}.  The exact same list of examples was rediscovered/re-interpreted by Thurston in terms of conical flat structures on the $2$ sphere \cite{ThurstonShapes}, see also \cite{SchwartzNotesonShapes}. There is an additional approach via algebraic geometry suggested by Hirzebruch and developed by him in collaboration with Barthels and H\"{o}fer \cite{BHH}. More examples have been discovered recently by Couwenberg, Heckman, and Looijenga using the Hirzebruch style techniques and by Deraux, Parker and Paupert using complex reflection group techniques \cite{CHL, DPP, DPP1, Deraux3d}.  But as of this writing there are only $22$ commensurability classes of non-arithmetic lattices known in $\SU(2,1)$ and only $2$ known in $\SU(3,1)$.   An obvious refinement of Question \ref{qtn:nonarithmeticSU} is the following:

\begin{qtn}
\label{qtn:nonarithmeticSUcomm}
For what values of $n$ do there exist infinitely many commensurablity classes of non-arithmetic lattice in $\SU(n,1)$?
\end{qtn}

\noindent We remark here that the approach via conical flat structures was extended by Veech and studied further by Ghazouani and Pirio \cite{VeechFlatSurf, GhazouaniPirio1}.  Regretably this approach does not yield more non-arithmetic examples.  It seems that the reach of this approach might be extended to be roughly equivalent to the reach of the approach via monodromy of hypergeometric functions, see \cite{GP2}.  There appears to be some consensus among experts is that the answer to both Question \ref{qtn:nonarithmeticSU} and Question \ref{qtn:nonarithmeticSUcomm} should be ``for all $n$", see e.g. \cite[Conjecture 10.8]{Kapovich}.  We point to a recent result of Esnault and Groechenig that indicates that complex hyperbolic lattices are in fact much more constrained than their real hyperbolic analogues \cite{EG1}.

\begin{thm}
  \label{thm:esnaultgroechenig}
Let $\Gamma$ be a lattice in $G=SU(n,1)$ for $n>1$.  Then $\Gamma$ is integral.  I.e. $\Gamma <G(k)$ for some number field $k$ and
if $\nu$ is any finite place of $k$ then $\Gamma < G(k_v)$ is pre-compact.
\end{thm}

The earliest non-arithmetic lattices in $\SO(n,1)$ for $n>2$ were constructed by Makarov and Vinberg by reflection group methods \cite{Makarov, VinbergFirst}.  It is known  by work of Vinberg that these methods will only produce non-arithmetic lattices in dimension less than $30$ \cite{Vinberg}.  The largest known non-arithmetic lattice produced by these methods is in dimension $18$ by Vinberg and the full limits of reflection group constructions is not well understood \cite{VinbergHigh}. We refer the reader to \cite{Belolipetsky} for a detailed survey. The following question seems natural:

\begin{qtn}
In what dimensions do there exist lattices in $\SO(n,1)$ or $\SU(n,1)$ that are commensurable to non-arithmetic reflection groups?  In what dimensions do there exist lattices in $\SO(n,1)$ or $\SU(n,1)$ that are commensurable to arithmetic reflection groups?
\end{qtn}

For the real hyperbolic setting, there are known upper bounds of $30$ for non-arithmetic lattices and $997$ for any lattices. The upper bound of $30$ also applies for arithmetic uniform hyperbolic lattices \cite{Vinberg, Belolipetsky}.  In the complex hyperbolic setting, there seem to be no known upper bounds, but a similar question recently appeared in e.g. \cite[Question 10.10]{Kapovich}.  For a much more detailed survey of reflection groups in hyperbolic spaces, see \cite{Belolipetsky}.

A dramatic result of Gromov and Piatetski-Shapiro vastly increased our stock of non-arithmetic lattices in $\SO(n,1)$
by an entirely new technique:

\begin{thm}[Gromov and Piatetski-Shapiro]
For each $n$ there exist infinitely many commensurability classes of non-arithmetic uniform and non-uniform lattices in $\SO(n,1)$.
\end{thm}

The construction in \cite{GPS} involves building hybrids of two arithmetic manifolds by cutting and pasting along totally geodesic codimension one submanifolds.  The key observation is that non-commensurable arithmetic manifolds can contain isometric totally geodesic codimension one submanifolds. This method has been extended and explored by many authors for a variety of purposes, see for example \cite{Agol, BelThom, SevenSamuraiShort, GelanderLevit}.  It has also been proposed that one might build non-arithmetic complex hyperbolic lattices using a variant of this method, though that proposal has largely been stymied by the lack of codimension one totally geodesic codimension one submanifolds.
The absence of codimension $1$ submanifolds makes it difficult to show that attempted ``hybrid" constructions yield discrete groups.  For more information see e.g. \cite{Paupert, PaupertWells, Wells} and
\cite[Conjecture 10.9]{Kapovich}.  We point out here that the results of Esnault and Groechenig discussed above implies that the ``inbreeding" variant of Agol and Belolipetsky-Thomson \cite{Agol, BelThom} cannot be adapted to produce non-arithmetic manifolds in the complex hyperbolic setting even if the original method of Gromov and Piatetski-Shapiro can be.  To me personally, this seems a very strong negative indication on the possibility of such an adaptation.  While the key observation of Agol
was made long after the paper of Gromov and Piateksi-Shapiro, the constructions are very similar.

In \cite{GPS}, Gromov and Piatetski-Shapiro ask the following intriguing question:

\begin{qtn}
\label{qtn:pieces}
Is it true that, in high enough dimensions, all lattices in $\SO(n,1)$ are built from sub-arithmetic pieces?
\end{qtn}

\noindent
The question is somewhat vague and \cite{GPS} also contains a group theoretic variant that is easily seen to be false for the examples of Agol
and Belolipetsky-Thomson, but a more precise starting point is:

\begin{qtn}
\label{question:tg}
For $n>3$, is it true that any non-arithmetic lattice in $\Gamma< \SO(n,1)$ intersects some conjugate of $\SO(n-1,1)$ in a lattice?
\end{qtn}

\noindent This is equivalent asking if every finite volume non-arithmetic hyperbolic manifold in dimension at least $4$ contains a closed codimension one totally geodesic submanifold.   Both reflection group constructions and all known variants of hybrid constructions contain such submanifolds.   It seems the consensus in the field is that the answer to this question should be no, but I know of no solid evidence for that belief.  In particular, starting in dimension $4$, Wang finiteness shows that the quantitative structure of hyperbolic manifolds is very different than in dimension $3$  \cite{Wang}.  And recent work of Gelander-Levit makes it at least plausible that variants of the hybrid and inbreeding constructions construct ``enough'' hyperbolic manifolds to capture all examples \cite{GelanderLevit}.  This is very different than the situation in dimension $3$ where hyperbolic Dehn surgery constructs many ``more'' hyperbolic manifolds and concretely exhibits the failure of Wang finiteness.  One can think of Questions \ref{qtn:pieces} and \ref{question:tg} as asking for particular qualitative reasons behind this difference in quantitative behavior.

 In the next subsection, I will discuss some approaches to giving a positive answer to this question or perhaps criterion for a positive answer for some examples building on ideas in \cite{BFMS}.  It is also not known to what extent the hybrid constructions and reflection group constructions build distinct examples. Some first results, indicating that the classes are different, are contained in \cite[Theorem 1.7]{FLMS} and in \cite[Theorem 1.5]{Mila}.

Given Theorem \ref{thm:bfms1}, it is reasonable to ask more detailed questions about the finite collection of totally geodesic submanifolds in a non-arithmetic
hyperbolic manifold.  A very reasonable question, on which first results have been obtained by Lindenstrauss and Mohammadi, is whether there is any bound
on the finite number in terms of geometric invariants of the hyperbolic manifold.  Their theorem, stated in Mohammadi's contribution to these proceedings as \cite[Theorem 7.4]{AmirICM}, gives
a result in the class of constructions of the style of Gromov and Piatetski-Shapiro in the case of three dimensional hyperbolic manifolds. A key ingredient
is the Angle Rigidity Theorem found in joint work of the author with Lafont, Miller and Stover \cite{FLMS}.  While the proof of finiteness in \cite{FLMS} is definitely superseded
by the one in \cite{BFMS}, it is more useful for proving bounds because it gives an explicit open set of ``impossible configurations" for the finite set of maximal
totally gedoesic submanifolds in many cut and paste constructions of hyperbolic manifolds. Roughly speaking, the key result in \cite{FLMS} says given a non-arithmetic manifold built in the manner of Gromov and Piatetski-Shapiro, any closed totally geodesic submanifold intersecting the ``cut and paste'' hypersurface must do so at a right angle. All other angles of intersection are forbidden.  It seems highly unlikely that these are the only ``impossible configurations".

\begin{qtn}
\label{qtn:configurations}
Find other restrictions on the posssible configurations of totally geodesic submanifolds in a non-arithmetic hyperbolic manifold.
\end{qtn}

It is worth mentioning that our understanding of lattices in $\SO(2,1)$ and $\SO(3,1)$ is both more developed and very different.
Lattices in $\SO(2,1)$ are completely classified, but there are many of them, with the typical isomorphism class of lattices
having many non-conjugate realizations as lattices, parameterized by moduli space.  In $\SO(3,1)$, Mostow rigidity means there
are no moduli spaces.  But Thurston-Jorgensen hyperbolic Dehn surgery still allows one to construct many ``more" examples of lattices,
including many that yield a negative answer to Question \ref{question:tg}.  There remains an interesting sense in which the answer
to Question \ref{qtn:pieces} could still be yes even for dimension $3$.

\begin{qtn}
\label{qtn:dehnsurgery}
Can every finite volume hyperbolic $3$-manifold be obtained as Dehn surgery on an arithmetic manifold?
\end{qtn}

\noindent To clarify the question, it is known that every finite volume hyperbolic $3$-manifold is obtained as a topological manifold  by Dehn surgery on some cover of the figure $8$ knot complement, which is known to be the only arithmetic knot complement \cite{HLM, ReidKnot}.  What is not known is whether one can obtain the geometric structure on the resulting three manifold as geometric deformation of the complete geometric structure on the arithmetic manifold on which one performs Dehn surgery.

We end this section by discussing an additional question that illustrates the difference in our knowledge of real hyperbolic manifolds in dimension $3$ and dimensions at least $4$ and then discuss an intriguing variant in complex hyperbolic geometry.  We call a group $\Gamma$ \emph{virtually large} if it has a finite index subgroup that surjects onto a non-abelian free group.  The following theorem of Lubotzky shows a strong connection between totally geodesic submanifolds and this property \cite{Lubotzky}.

\begin{thm}\label{thm:lubotzky}
Let $M$ be a finite volume hyperbolic manifold.  Then if $M$ admits a closed codimension $1$ totally geodesic submanifold, the fundamental group of $M$ is virtually large.
\end{thm}

\noindent It follows from work of Agol that in dimension $3$ all finite volume hyperbolic manifolds have virtually large fundamental group \cite{Agol2}.  But starting
in dimension $5$ there are explicit examples where we do not know if the fundamental group is virtually large.  These are the so-called \emph{second type} arithmetic groups
constructed using quaternion algebras.  In fact, effectively the only way we know to show a hyperbolic manifold of dimension $4$ or higher has virtually large fundamental group
is to apply Lubotzky's theorem.  We ask the following question in an intentionaly provocative manner:

\begin{qtn}
\label{qtn:large}
If $M$ is a finite volume hyperbolic manifold of dimension at least $4$ is having virtually large fundamental group equivalent to having a closed totally geodesic submanifold of codimension $1$?  \end{qtn}

\noindent There seems to be something close to a consensus that the answer to this question should be ``no'' in general and all lattices in $\SO(n,1)$ should be virtually large.  However, I know of no strong evidence for this belief other than the results in dimension $3$. There is even one potential strategy for proving a negative answer to Question \ref{qtn:large}, arising from Agol's work on the virtual Haken conjecture \cite{AgolFiber}.

\begin{qtn}
\label{qtn:cubes}
If $M$ is a finite volume hyperbolic manifold, can we cubuluate $M$?
\end{qtn}

The answer is yes in dimension $3$ by work of Kahn-Markovic \cite{KahnMarkovic}.  Agol's work then implies that this cubulation is special, which is more than
enough to prove largeness.  The concrete challenge in this setting is to prove that \emph{second type} arithmetic hyperbolic manifolds can be cubulated.  That the first
type can be is shown by Bergeron, Haglund and Wise in \cite{MR2776645}.  It seems well known to experts that all inbreeding and hybrid examples can also be cubulated
but there does not appear to be an explicit statement in the literature.  However all of these approaches to cubulating hyperbolic manifolds of dimension at least $4$ depend on the existence of totally geodesic submanifolds of codimension $1$.  To cubulate in cases where codimension $1$ totally geodesic submanifolds do not exist, work in dimension $3$ suggests one should find some other convex or quasi-convex submanifolds or subspaces of codimension $1$.  Due to density of the commensurator, in the arithmetic setting, it suffices to one such manifold.

It is interesting in the context of this article to also summarize what is known about largeness for lattices in $\SU(n,1)$.  In this case all
known constructions don't produce largeness directly.  Instead one constructs a holomorphic retract onto a totally geodesic Riemann subsurface as a
forgotful map in terms of the Deligne-Mostow hypergeometric monodromy construction \cite{DerauxForgetful}.  It is possible that only these particular Deligne-Mostow constructions yield virtually large lattices in $\SU(n,1)$, but here we only ask a weaker question.

\begin{qtn}
\label{qtn:forget}
Let $M$ be a finite volume complex hyperbolic manifold.  If the fundamental group of $M$ is large, does $M$ admit a holomorphic retract to totally geodesic Riemann subsurface?
\end{qtn}

We remark here that by a result of Delzant and Py, complex hyperbolic manifolds of (complex) dimension at least $2$ are not cubulated, so the approach
to a negative answer to Question \ref{qtn:large} is not available in this context \cite{DelzantPy}.

\subsection{Dynamics and (non)-arithmeticity}

The work in \cite{BFMS, BFMS2} gives a broader set of tools for determining when a manifold is arithmetic.  We discuss some further repercussions of these ideas here, though they have less decisive corollaries than the main results of those papers.  We restrict here to the real hyperbolic setting for simplicity, so throughout this section $G=\SO(n,1)$ with $n>2$.  We will write $W=SO(n-1,1)$, though for some technical results below, there are analogous statements when $W=SO(k,1)$ for any $1<k<n-1$.

It is known that given a lattice $\Gamma<G$, there exists a number field $\ell$ and an $\ell$ structure on $G$ such that $\Gamma< G(\ell)$.  We note that there is a collection of valuations on $\ell$ and we write $\nu_0$ for the valuation for which $G(\ell_{\nu_0})=\SO(n,1)$ in such a way that we get the given lattice embedding of $\Gamma$ into $\SO(n,1)$.

\begin{defn}
\label{defn:arithobstruction}
Given a valuation $\nu$ on $\ell$ not equivalent to $\nu_0$, the \emph{arithmeticity obstruction} defined by $\nu$ is the embedding
$\rho_{\nu}:\Gamma \rightarrow G(\ell_{\nu})=H$.  We say the arithmeticity obstruction
vanishes if $\Gamma$ is pre-compact in $H$.
\end{defn}

\noindent It was observed by Margulis that if all arithmeticity obstructions vanish, then $\Gamma$ is arithmetic.  If all non-archimedean arithmeticity obstructions vanish,
then $\Gamma$  is called \emph{integral}.  Note that $\Gamma$ is automatically Zariski dense in any $G(\ell_{\nu})$.  The standard and surprising technique for showing that arithmeticity obstructions vanish is to assume they don't and use a superrigidity theorem to see that this implies that $\rho_{\nu}$ extends to a continuous representation from $G$ to $H$. This is easily seen to be impossible.

We now state the superrigidity theorem from \cite{BFMS} that is used there to show arithmeticity obstructions vanish. We restrict attention to the case of $\SO(n,1)$ for simplicity.   This theorem is stated in terms of the existence of a certain $W$ invariant measure on a certain flat bundle over $G/\Gamma$.  The goal of this section is to explain that it is possible to use other dynamical results to weaken that hypothesis.  Finding optimal hypothesis and more applications seems potentially important to understanding the questions discussed in the last subsection.

In this theorem we consider a local field $k$ and a $k$-algebraic group $\bfH$ satisfying one additional condition. Let $P$ be a minimal parabolic subgroup of $G$ and $U$ its unipotent radical. A pair consisting of a local field $k$ and a $k$-algebraic group $\bfH$ is said to be \emph{compatible} with $G$ if for every nontrivial $k$-subgroup $\bfJ<\bfH$ and any continuous homomorphism $\tau:P \rightarrow N_{\bfH}(\bfJ)/\bfJ(k)$, where $N_{\bfH}(\bfJ)$ is the normalizer of $\bfJ$ in $\bfH$, we have that the Zariski closure of $\tau(U')$ coincides with the Zariski closure of $\tau(U)$ for every nontrivial subgroup $U'<U$.  That this condition is satisfied by any group arising in an arithmeticity obstruction for a lattice in $\SO(n,1)$ is a key point in \cite{BFMS}.  The fact this this is already no longer true for $\SU(n,1)$ introduces new difficulties overcome in \cite{BFMS2} that we do not discuss here.

\begin{thm}
\label{theorem:superrigiditydichotomy}
Let $G$ be $\SO_0(n, 1)$ for $n \ge 3$, $W<G$ be a noncompact simple subgroup, and $\Gamma < G$ be a lattice. Suppose that $k$ is a local field and $\bfH$ is a connected $k$-algebraic group such that the pair consisting of $k$ and $\bfH$ is compatible with $G$. Finally, let $\rho: \Gamma \rightarrow \bfH(k)$ be a homomorphism with unbounded, Zariski dense image. If there exist a $k$-rational faithful irreducible representation $\bfH \to \SL(V)$ on a $k$-vector space $V$ and a $W$-invariant measure $\nu$ on $(G\times \mathbb{P}(V))/\Gamma$ that projects to Haar measure on $G / \Gam$, then $\rho$ extends to a continuous homomorphism from $G$ to $\bfH(k)$.
\end{thm}

We retain the assumptions made on $k$ and $\bfH$ and discuss how one can weaken the assumption in Theorem \ref{theorem:superrigiditydichotomy} using results of
Eskin-Bonatti-Wilkinson.  We first motivate the connection by stating a variant that follows easily.

\begin{thm}[BFMS reformulated]
\label{thm:bfmsprime}
Take the assumptions of Theorem \ref{theorem:superrigiditydichotomy} and replace the assumption of the existence of the invariant measure with the assumption that
the $W$ action on $(G \times V)/\Gamma$ is not irreducible.  Then $\rho$ extends to a continuous homomorphism from $G$ to $\bfH(k)$.
\end{thm}

The two statements are equivalent by standard techniques.   The point is that either the $W$ invariant measure or the $W$ invariant subspace are used to produce
a $\Gamma$-equivariant measurable map $\phi: W \backslash G \rightarrow (\bfH/\bfL)(k)$ for some proper algebraic subgroup $\bfL <\bfH$. The existence of this map is in fact
equivalent to the existence of either such a subpace for some choice of $V$ or such a measure for some (a priori different) choice of $V$.  In this language, the first step of the work in \cite{BFMS} is using the existence of infinitely many maximal $W$ orbit closures to produce an invariant measure for the $W$ action on $(G\times \mathbb{P}(V))/\Gamma$ for some well chosen $V$.  It is worth pointing out that there are other criteria for irreduciblity in the literature, here we focus on one due to Bonatti, Eskin and Wilkinson \cite{BEW}.  To obtain a simple statement, we let $A<P$ be the Cartan subgroup and choose the $H$ representation on $V$ such that the first Lyapunov exponent for $A$ on $V$-bundle $(G \times V)/\Gamma$ is simple.  It is easy to verify that one can choose such a $V$ using tensor constructions.   We let $P<G$ be the parabolic and $P_W= P \cap W$ the parabolic in $W$. Combining the main result of \cite{BEW} with those of \cite{BFMS} we obtain the following:

\begin{thm}\label{thm:technical}
With the common assumptions of the last two theorems if there is more than one $P_W$ invariant measure on $(G \times \bfP(V))/\Gamma$ projecting to Haar measure on $G/\Gamma$,
then  $\rho$ extends to a continuous homomorpism from $G$ to $\bfH(k)$.
\end{thm}

\noindent This is a straightforward concatenation of \cite[Theorem 1.1]{BEW} and Theorem \ref{thm:bfmsprime}.  Observe first that the $V$-bundle
$(G \times V)/\Gamma$ is $G$ irreducible by hypothesis.  The assumption of the second invariant measure in Theorem \ref{thm:technical} and \cite[Theorem 1.1]{BEW}
imply that this bundle is not irreducible for $W$. We then apply Theorem \ref{thm:bfmsprime}.

Let $E_1(x)$ be the subspace of $V$ corresponding to the first Lyapunov exponent for $A$ and write $(G \times E_1)/\Gamma$ for the corresponding subbundle of
$(G \times V)/\Gamma$.  It is observed in \cite{BEW} that $(G \times E_1)/\Gamma$ is $P$ invariant.  Using that $E_1$ is one dimensional, this yields a $P$-invariant section
$s$ of the bundle $(G \times \bfP(V))/\Gamma$ and we can push the Haar measure on $G/\Gamma$ forward along this section to build a $P$ invariant measure on $(G\times \bfP(V))/\Gamma$.
We note that this measure cannot be $W$ invariant since $P$ and $W$ together generate $G$ and it is easy to see that there is no $G$ invariant measure on $(G\times \bfP(V))/\Gamma$.
Since $P_W<P$, This shows the existence of a $W$ invariant measure on $(G\times \bfP(V))/\Gamma$ projecting to Haar implies the existence of at least two $P_W$ invariant measures on $(G\times \bfP(V))/\Gamma$ projecting to Haar.  It seems a priori easier to produce $P_W$ invariant measures since $P_W$ is amenable and therefore one can average.  Producing a measure that is not also $P$ invariant becomes the challenge.

We also note here that this entire discussion on varying hypotheses, applies in any case where $G$ is a rank one group and $W<G$ is a simple subgroup, we only require $G=\SO(n,1)$ to use the superrigidity theorem from \cite{BFMS} to extend the representation. So the discussion adapts easily to the case of $\SU(n,1)$ by the results
of \cite{BFMS2}.   There is one additional fact that is special to $\SO(n,1)$ and relates to Questions \ref{question:tg} and \ref{qtn:large}, namely that if the hyperbolic manifold $K\backslash G/\Gamma$ has no totally geodesic manifolds of codimension one and $W=\SO(n-1,1)$ then
the $P_W$ action on $G/\Gamma$ is uniquely ergodic.  So in this setting we have the following.

\begin{cor}
\label{cor:uniqueergodic}
Adding the assumption on totally geodesic submanifolds above to Theorem \ref{thm:technical}, if $P_W$ is not uniquely ergodic on $(G\times \bfP(V))/\Gamma$ then $\rho$ extends.
\end{cor}

It is somewhat surprising that arithmeticity is in this context associated with the existence of additional ergodic measures.  Since $\Gamma$ is finitely generated, it is easy to see that only finitely many valuations of $\ell$ are relevant to arithmeticity of $\Gamma$.  So arithmeticity of $\Gamma$ follows from the failure of unique ergodicity for $P_W$ in finitely many dynamical systems for $G$.  We note that the paper of Bonatti-Eskin-Wilkinson gives somewhat more explicit information than just about the number of measures.  There is one "obvious" ergodic measure supported on the line in $P(V)$ corresponding to the first Lyapunov exponent.  To force vanishing of the arithmeticity obstruction simply requires finding any $P_W$ invariant measure where disintegration along fibers is not supported on that line.

In the context of Question \ref{qtn:large}, the possibilities for applying \ref{cor:uniqueergodic} are much broader.  Given a lattice $\Gamma$ that is large, one has a homomorphism $\Gamma \rightarrow F_2$ which one can compose with any (irreducible) representation of $F_2$ on a vector space $V$ to obtain a $G$ space $(G \times \bfP(V))/\Gamma$
to which one can attempt to apply Corollary \ref{cor:uniqueergodic}. Since a representation factoring through $F_2$ cannot extend, the corollary actually asserts that largeness implies unique ergodicity of $P_W$ in many dynamical systems. Of course, all lattices in  $PSL(2,\mathbb{C})$ are known to be virtually large and there are many where the locally symmetric space has no closed totally geodesic surfaces, so for that choice of $G$ all of these constructions build many dynamical systems with a unique invariant measures for the corresponding $P_W$.  In this case $P_W$ is just the $ax+b$ group.  It is not at this moment clear what structure one would like to exploit to show a difference between dimesnion $3$ and higher dimensions.

One also might attempt to use Theorem \ref{thm:technical} to either reprove the arithmeticity of lattices in $SP(n,1)$ and $F_4^{-20}$ by ergodic theoretic methods or to reprove Theorem \ref{thm:esnaultgroechenig} by ergodic theoretic methods.  In addition, one could attempt a proof of Mostow rigidity using this circle of ideas.  Certainly the ideas are close to those that allow Margulis to prove Mostow rigidity as a consequence of superrigidity in higher rank.

We note that while Margulis and Mohammadi suggested an alternate approach to the theorem in \cite{BFMS} and showed their approach could succeed in the case where $G=\SO(3,1)$, their approach does not yield any of the observations in this subsection but only produces results in the presence of infinitely many maximal closed totally geodesic submanifolds.


\begin{ack}
I am happy to thank all of my collaborators, without whom the works discussed here would not exist.  Particular thanks to Nick Miller for a detailed reading of an
earlier draft and to Ian Agol and Matthew Stover for interesting conversations about largeness.
\end{ack}

\begin{funding}
This work was partially supported by a Fellowship from the Simons Foundation, by a visiting professorship at the Miller Research Institute and by NSF DMS-1906107.
\end{funding}


\bibliographystyle{emss}
\bibliography{bibliography}









\end{document}